\documentclass[a4paper, 10pt, conference]{ieeeconf}      

\IEEEoverridecommandlockouts                              

\usepackage{graphics} 
\usepackage{epsfig} 
\usepackage{mathptmx} 
\usepackage{times} 
\usepackage{amsmath} 
\usepackage{amssymb}  
\usepackage{xcolor}         
\usepackage{subcaption}
\usepackage{floatrow}
\usepackage{placeins}
\usepackage{blindtext}
\usepackage{algorithm}
\usepackage{algpseudocode}
\usepackage{multirow}
\usepackage{wrapfig}
\usepackage{hyperref}
\usepackage{booktabs}

\title{\LARGE \bf Using Laplace Transform To Optimize the Hallucination of Generation Models}

\author{\title{\LARGE \bf Using Laplace Transform To Optimize the Hallucination of Generation Models}}

\author{
Cheng Kang$^{1}$, Xinye Chen$^{2}$, Daniel Novak$^{1}$, and Xujing Yao$^{3}$\textsuperscript{*}\thanks{*Corresponding author: Xujing Yao (xjyao@njtech.edu.cn).}%
\thanks{This work was supported by the inEXASCALE (101075632), Czech Technical University in Prague (SGS22/165/OHK3/3T/13), and the Research Centre for Informatics (CZ.02.1.01/0.0/0.0/160 19/0000765).}%
\thanks{$^{1}$Cheng Kang and Daniel Novak are with the Department of Cybernetics, Czech Technical University in Prague, 16000 Praha 6, Prague, Czech Republic {\tt\small \{kangchen,xnovakd1\}@fel.cvut.cz}}%
\thanks{$^{2}$Xinye Chen is with Department of Numerical Mathematics, Charles University, 11000 Praha 1, Prague, Czech Republic {\tt\small xinye.chen@mff.cuni.cz}}%
\thanks{$^{3}$Xujing Yao is with the College of Computer and Information Engineering, Nanjing Tech University, Nanjing, Jiangsu 211816, China {\tt\small xjyao@njtech.edu.cn}}%
}

\begin{document}

\maketitle
\thispagestyle{empty}
\pagestyle{empty}

\begin{abstract}
To explore the feasibility of avoiding the confident error (or hallucination) of generation models (GMs), we formalise the system of GMs as a class of stochastic dynamical systems through the lens of control theory. Numerous factors can be attributed to the hallucination of the learning process of GMs, utilising knowledge of control theory allows us to analyse their system functions and system responses. Due to the high complexity of GMs when using various optimization methods, we cannot figure out their solution of Laplace transform, but from a macroscopic perspective, simulating the source response provides a virtual way to address the hallucination of GMs. We also find that the training progress is consistent with the corresponding system response, which offers us a useful way to develop a better optimization component. Finally, the hallucination problem of GMs is fundamentally optimized by using Laplace transform analysis. 
\end{abstract}

\section{INTRODUCTION}

Generation models, especially, large language models (LLMs) and large vision generation models (LVGMs), are suffering from the hallucination problem. Speaking of language generation and vision-language generation models, hallucination is referred to as a situation where the model generates content that is not based on factual or accurate information \cite{rawte2023survey}. Hallucination in vision generation models is a phenomenon that models generate samples that lie completely out of the training distribution of the model \cite{aithal2024understanding}. One conversational study investigated the original of hallucination, and they found that hallucination is a prevalent issue in both dialog benchmarks and models \cite{dziri2022origin}. Recent studies of diffusion generation models indicates that the hallucination phenomenon has a ''mode collapse'' \cite{zhang2018convergence} which leads to a loss in diversity in the sampled distribution.

So how does existing studies address hallucination? Focusing on the microscopic level, a large commonality in the majority of prior work (both in language generation models \cite{shuster2021retrieval,dziri2021neural} and vision generation models \cite{zhang2018convergence,kim2024tackling,leng2024mitigating}) seeks to address hallucination by ameliorating the model. Retrieval augmented generation (RAG) improved the reliability of LLMs by reading extra knowledge from new data \cite{shuster2021retrieval}. For language-vision models, existing studies \cite{hu2023ciem,zhai2023halle} mitigated the issue of hallucination in image captioning models with a predominant focus on the “presence of objects”, specifically whether or not the objects depicted in a given image are accurately described by the text that is generated by the mode. Most advanced generation adversarial models (GANs) assembled some quality components \cite{zhou2023gan} that wonderfully optimize the generation results. Paying attention to the local out-of-distribution (OOD) regions of images when using denoising diffusion probabilistic models (DDPM) to translate one image to another image,  multiple local diffusion processes proved to reduce the hallucination of generated images \cite{kim2024tackling}. However, to the best of our knowledge, apart from the problem of dataset, there is no research and discussion on the stability of these learning systems, as well as the fundamental reason why generation models have confident error or hallucination in general.

Can we fundamentally figure out why generation models have hallucination from a macroscopic perspective? Numerous advanced components, especially, optimization methods \cite{wang2020pid} and loss functions \cite{chu2017cyclegan}, have improved generation models, and their utilisation will affect the performance of GMs. Conducting a quantitative analysis on each of them can pave the way for their efficiency and development. In one node of the learning system, and in terms of analyzing a single component by computing their Laplace transform, their stability and convergence \cite{nise2020control} can be evaluated to finally ensure that GMs are able to generate a reliable and stable output.

How well does Laplace transform address the hallucination of GMs? Based on the Laplace transform analysis \cite{nise2020control}, we find that using a proper optimizer (such as, using Adam optimizer on GAN \cite{2014Generative} and DDPM \cite{ho2020denoising}, using FuzzyPID optimizer \cite{yao2024automatic} on CycleGAN \cite{chu2017cyclegan}), can benefit the learning process of GMs. The use of a proper component can ensure a harmony of the whole system, which can visually tell us which component should be used in which learning system. In this study, to guarantee the stability and convergence of GMs, we find that there are two ways that can address the hallucination of GMs via Laplace transform: \textbf{(1)} using Laplace transform to develop a proper component; \textbf{(2)} using Laplace transform to design a better learning system. 

Contributions are shown below.

\begin{itemize}

\item Laplace transform analysis finds the fundamental reason why generation models have hallucination.
    
\item Optimization component should  specifically designed to match the generation models;

\item The practical performance of every optimization component is highly consistent with their system response.

\end{itemize}

\section{The Connection Between Control Theory and Learning Progress of GMs}

\begin{figure}
    \centering 
    \includegraphics[scale=0.28]{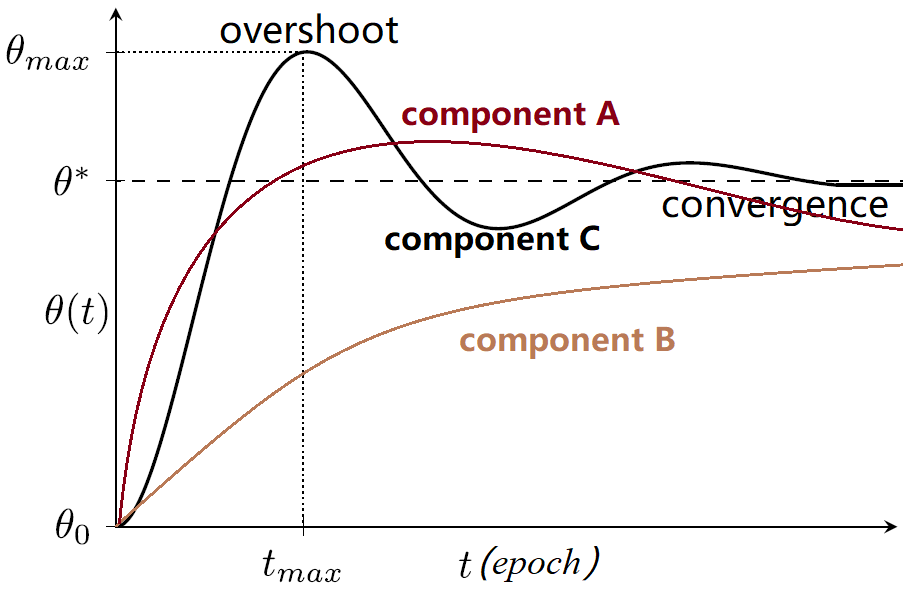} 
    \caption{Evolution process of the weight $\theta({t})$ for each component under the generation task.} 
    \label{Figure-theory} 
\end{figure}

In Figure \ref{Figure-theory}, $\theta_{max}$ is the maximum overshoot of training model, and $t_{max}$ is the needed epoch or step when the training model reaches the strongest vibration. We assume that the overshoot and vibration is the source of hallucination. The optimal condition is that the training model has a proper maximum vibration (no vibration will make model learn slowly, for example, component B in Figure \ref{Figure-theory}), a shorter climbing time $t_{max}$, a visual robustness, and a reliable and convergent result (the learning process will result in a long time to convergence or no convergence, for example, component A in Figure \ref{Figure-theory}). Compared to the real GM itself, so as to obtain the optimal $\theta^{\ast}$, we should guarantee that the learning progress is associated with a stable training progress for a better result with a shorter training time. For a clear simulation analysis, we assume that the system function of every GM is a two-order system, and their system function can be given by $1/(s^2+P_{a}s+P_{b})$), where $P$, $P_{a}$ and $P_{b}$ are the corresponding variables, respectively.

\section{The Laplace Transform of Optimization Methods}

In this section, we review several widely used optimizers, such as SGD \cite{robbins1951stochastic,cotter2011better,zhou2017convergence}, SGDM \cite{qian1999momentum,liu2020improved}, Adam \cite{kingma2014Adam,bock2018improvement}, PID-optimizer \cite{wang2020pid} and Gaussian LPF-SGD \cite{bisla2022low}.

\subsection{SGD Optimization Method}

The parameter update rule of SGD from iteration $t$ to $t+1$ is determined by
\begin{equation} \label{Eq21}
\theta_{t+1}=\theta_t-r \frac{\partial L_t}{\partial \theta_t}, 
\end{equation}
where $r$ is the learning rate. We now regard the gradient $\partial L_t / \partial \theta_t$ as error $e(t)$ in the PID control system \cite{wang2020pid}. Compared to the PID controller, we find that SGD can be viewed as one type of $P$ controller with $K_p=r$. The system function of SGD finally becomes:
\begin{equation} \label{Eq22}
\theta_{SGD}(s)=r.
\end{equation}

\subsection{SGDM Optimization Method}

By leveraging historical gradients, SGDM trains a ANN more swiftly than SGD does. The rule of SGDM updating parameter is given by
\begin{equation} \label{Eq23}
\left\{\begin{array}{l}
V_{t+1}=\alpha V_t-r \partial L_t / \partial \theta_t \\
\theta_{t+1}=\theta_t+V_{t+1}
\end{array}\right.
\end{equation}
where $V_t$ is a term that accumulates historical gradients. $\alpha \in(0,1)$ is the factor that balances the past and current gradients. It is usually set to $0.9$. Dividing two sides of the Equation \ref{Eq23} by $\alpha^{t+1}$, we get: 
\begin{equation} \label{Eq24}
\frac{V_{t+1}}{\alpha^{t+1}}=\frac{V_t}{\alpha^t}-\frac{r}{\alpha^{t+1}} \frac{\partial L_t }{\partial \theta_t} .
\end{equation}

\noindent Finally, we get $\theta_{t+1}$ as follow by iteration:
\begin{equation} \label{Eq25}
\theta_{t+1}-\theta_t=-r \frac{\partial L_t}{\partial \theta_t}-r \sum_{i=0}^{t-1} \alpha^{t-i} \frac{\partial L_i}{\partial \theta_i}.
\end{equation}

\noindent SGDM actually is a PI controller with $K_p=r$ and $K_i=r\alpha^{t-i}$. The system function of SGDM is:
\begin{equation} \label{Eq26}
\theta_{SGDM}(s)=r + \frac{r}{s} \cdot \frac{1}{s-ln(\alpha)}.
\end{equation}

\subsection{PID Optimization Method}

SGD and SGDM can be respectively viewed as P and PI controller \cite{wang2020pid}. Given that training is often conducted in a mini-batch manner, the learning process is very easy to introduce noise when computing gradients. The proposed PID optimizer \cite{wang2020pid} updates network parameter $\theta$ in iteration $(t+1)$ by
\begin{equation} \label{Eq27}
\left\{\begin{array}{l}
V_{t+1}=\alpha V_t-r \partial L_t / \partial \theta_t \\
D_{t+1}=\alpha D_t+(1-\alpha)\left(\partial L_t / \partial \theta_t-\partial L_{t-1} / \partial \theta_{t-1}\right) \\
\theta_{t+1}=\theta_t+V_{t+1}+K_d D_{t+1} .
\end{array}\right.
\end{equation}

\noindent Thus, the $\theta_{t+1}$ using PID optimizer is described as follow by iteration:
\begin{equation} \label{Eq28}
\theta_{t+1}-\theta_t=-r \frac{\partial L_t}{\partial \theta_t} - r \sum_{i=0}^{t-1} \alpha^{t-i} \frac{\partial L_i}{\partial \theta_i} - r K_{d} \left( \frac{\partial L_i}{\partial \theta_i} - \frac{\partial L_{i-1}}{\partial \theta_{i-1}} \right).
\end{equation}

\noindent where $K_{d} \left( \frac{\partial L_i}{\partial \theta_i} - \frac{\partial L_{i-1}}{\partial \theta_{i-1}} \right)$ is the D component of the PID controller. The system function of PID is:

\begin{equation} \label{Eq29}
\theta_{PID}(s) = r + \frac{r}{s} \cdot \frac{1}{s-ln(\alpha)} + K_{d} s.
\end{equation}

When setting the hyperparameter $\alpha=1.0$, we can get the vanilla PID optimizer: $K_{p}=r$, $K_{i}=r$ and $K_{d}=r \cdot K_{d}$

\subsection{Adam Optimization Method}

Based on adaptive estimates of lower-order moments, Adam algorithm adaptively adjusts the stochastic gradients, and it can be summarized as below:

\begin{equation} \label{Eq30}
\left\{\begin{array}{l}
m_{t+1}=\beta_{1} m_t + \left(1-\beta_{1}\right) \partial L_t / \partial \theta_t \\
v_{t+1}=\beta_{2} v_t + \left(1-\beta_{1}\right) \partial L_t / \partial \theta_t \\
\hat{m}_{t+1}=m_{t} / \left( 1 - \beta^{t}_{1} \right) \\
\hat{v}_{t+1}=v_{t} / \left( 1 - \beta^{t}_{2} \right) \\
\theta_{t+1}=\theta_t + \alpha \hat{m}_{t} / \left( \sqrt{\hat{v}_{t+1}} + \epsilon \right)
\end{array}\right.
\end{equation}
where $m_{t}$ is the first moment estimate at timestep $t$, and $v_{t}$ is the second raw moment estimate. The default sets of learning rate $\alpha$, hyperparameters $\beta_{1}$, $\beta_{2}$ and $\epsilon$ are respectively $0.001$, $0.9$, $0.999$ and $10^{ - 8}$.

The iteration of $\theta_{t+1}$ using the Adam optimizer is described as follows:

\begin{equation} \label{Eq31}
\begin{aligned}
\theta_{t+1} & = \theta_{t} - r \cdot \frac{\widehat{m}_{t}}{\sqrt{\widehat{v}_{t}} + \epsilon} \\
& = \theta_{t} - r \cdot \frac{\frac{\sum_{i=0}^{t}\beta_{1}^{t-i}(\partial L_{i} / \partial \theta_{i})}{ \sum_{i=1}^{t}\beta_1^{i-1} }}{\sqrt{{\frac{\sum_{i=1}^{t}\beta_{2}^{t-i}(\partial L_{i} / \partial \theta_{i})^2}{ \sum_{i=1}^{t}\beta_2^{i-1} }}} + \epsilon} \\
& = \theta_{t} - r \cdot \frac{1}{M} \beta_{1}^{0} \frac{\partial L_{t}}{\partial \theta_{t}} - r \cdot \frac{1}{M} \sum_{i=0}^{t-1}\beta_{1}^{t-1-i} \frac{\partial L_{i}}{\partial \theta_{i}},
\end{aligned}
\end{equation}
where $M$ is the adaptive part of Adam, and its formula is:

\begin{equation} \label{Eq32}
M = \frac{1}{\sqrt{\frac{ \sum_{i=0}^{t}\beta _{2}^{t-i}( \partial L_t / \partial \theta_t)^{2} }{ \sum_{i=0}^{t}\beta_{2}^{i-1} }} + \epsilon } \cdot \frac{ 1 }{ \sum_{i=0}^{t} \beta_{1}^{i-1} }.
\end{equation}

Compared to Equation \ref{Eq25}, the adaptive component $M$ of Adam optimizer plays an critical role on adapting the learning system. We cannot derive the system function of Adam, as the high complexity of $M$. Finally, we directly use the same S function in Simulink and get its system response on above mentioned generation models.

\subsection{Filter Processed SGD Optimization Method}
Although Gaussian LPF-SGD outperforms other SGD variants, we are still not clear which part it has filtered, for example, high frequency, low frequency or any band frequency parts. In this study, we summarize the SGD learning process under the processing of filters as below:

\begin{equation} \label{Eq33}
\left\{\begin{array}{l}
\widehat{\frac{\partial L_t}{\partial \theta_t}} = \frac{\partial L_t}{\partial \theta_t} + \frac{\partial \left( \int_{-\infty}^{\infty} L(\theta_{t}-\tau) H(\tau) d \tau \right)}{\partial \theta_t} \\
\theta_{t+1} = \theta_t - r \widehat{\frac{\partial L_t}{\partial \theta_t}}
\end{array}\right.
\end{equation}
where $H$ is a Gaussian kernel and $L(\theta_{t})$ is the loss function of the training process in GLPF-SGD \cite{bisla2022low}. The $\theta_{t+1}$ using $Filter$ processed SGD optimizer is described as follow by iteration:

\begin{equation} \label{Eq34}
\begin{aligned}
\theta_{t+1} & = \theta_t - r \frac{\partial L_t}{\partial \theta_t} + r \frac{\partial \left( \int_{-\infty}^{\infty} L(\theta_{t}-\tau) H(\tau) d \tau \right)}{\partial \theta_t}\\
& = \theta_{t} - r \frac{\partial L_{t}}{\partial \theta_{t}} + r \frac{\partial \left( L {\circledast} H \right)}{\partial \theta_{t}}\\
& = \theta_t - r \frac{\partial L_t}{\partial \theta_t} + r \frac{1}{G}\frac{\sum_{i=0}^N \partial \left( L(\theta_t-\tau_i) H(\tau_i)\right)}{\partial \theta_t},
\end{aligned}
\end{equation}
where $G$ is the gain of the filter $H$ with the order of $N$, and $\circledast$ is the convolution process. Finally, the system function of filter processed SGD becomes:
\begin{equation} \label{Eq35}
\theta_{FP-SGD}(s) = r \left(G \cdot \frac{\prod_{i=0}^{m} \left(s+h_{i}\right)}{\prod_{j=0}^{n} \left(s+l_{j}\right)} \right).
\end{equation}


In this study, to analyse which frequency parts are beneficial to the training, we used a second-order Infinite Impulse Response (IIR) filter instead of the Gaussian kernel filter. By approximately setting the cutoff frequency at half, we imply a low-pass filter ranging from 0 Hz to half the sampling rate and a high-pass filter from half the sampling rate up to the sampling rate. Consequently, knowledge of the exact sampling rate is unnecessary, and essentially, it remains unobtainable.

\subsection{FuzzyPID Optimization Method}

Based on PID optimizer \cite{wang2020pid}, a designed PID controller which is optimized by fuzzy logic can make the training process more stable while keeping the dominant attribute of models. For instance, the ability of reaching the quick convergent result.

There are two key factors which affect the performance of the Fuzzy PID optimizer: (1) the selection of Fuzzy Universe Range $[-\varphi, \varphi]$ and (2) Membership Function Type $f_{m}$. 
\begin{gather} \label{eq10}
\begin{split}
& \widehat{K}_{\mathrm{P,I,D}}={K}_{\mathrm{P,I,D}}+\Delta K_{\mathrm{P,I,D}}
\end{split}
\end{gather}
\begin{gather} \label{eq11}
\begin{split}
& \Delta K_{\mathrm{P,I,D}} = D_{f}(E(s),Ec(s)) \cdot K_{\mathrm{P,I,D}}\\
& D_{f}(s) = f_{m}(round(-\varphi, \varphi, s)),
\end{split}
\end{gather}
where $D_{f}$ is the defuzzy process, $\Delta{K}_{\mathrm{P,I,D}}$ denotes the default gain coefficients of ${K}_{\mathrm{P}}$, ${K}_{\mathrm{I}}$ and ${K}_{\mathrm{D}}$ before the modification. $E(s)$ is referred to as the back-propagation error, and $Ec(s)$ is the difference between the $Laplace$ of $e(t)$ and $e(t-1)$. The Laplace function of this model $\theta(s)$ eventually becomes:
\begin{equation} \label{Eq12}
\theta(s) = \frac{\widehat{K}_{d}s^2+\widehat{K}_{p}s+\widehat{K}_{i}}{\widehat{K}_{d}s^2+(\widehat{K}_{p}+1)s+\widehat{K}_{i}} \cdot \frac{\theta^\ast}{s}.
\end{equation}
\noindent where $\widehat{K}_{p}$, $\widehat{K}_{i}$ and $\widehat{K}_{d}$ should be processed under the fuzzy logic. By carefully selecting the learning rate $r$, $\theta(s)$ becomes a stable system. The PID \cite{ang2005pid} and Fuzzy PID \cite{tang2001optimal} controllers have been used to control a feedback system by exploiting the present, past, and future information of prediction error. The advantages of a fuzzy PID controller includes that it can provide different response levels to non-linear variations in a system. At the same time, the fuzzy PID controller can function as well as a standard PID controller in a system where variation is predictable.

\section{The Laplace Transform of GAN}

\begin{figure}[ht]
\centering
\includegraphics[scale=0.44]{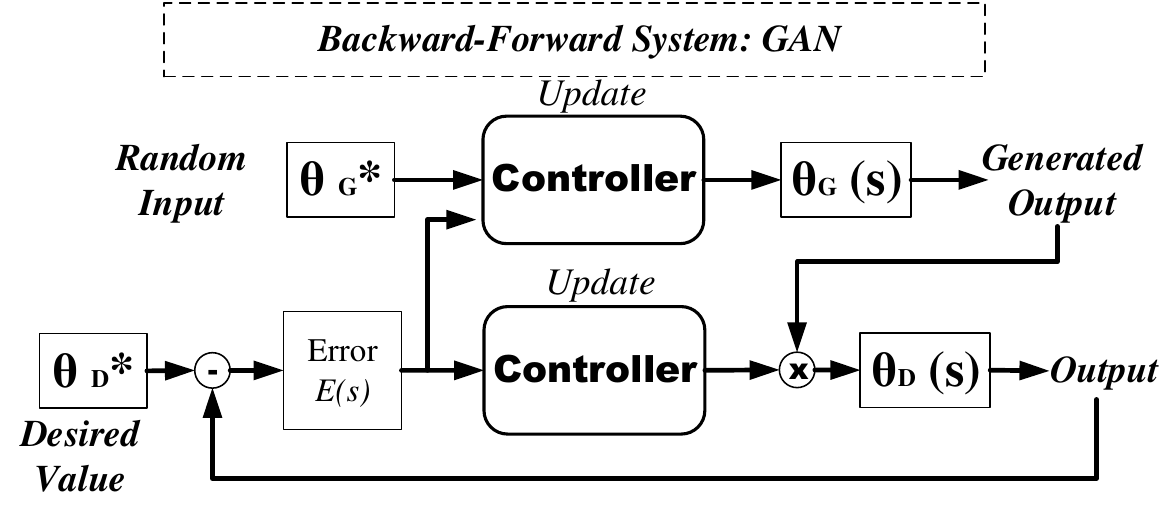}
\captionof{figure}{The control system of classical GAN.}
\label{Figure1}
\end{figure}

\begin{figure*}[htp!]
    \centering 
    \begin{subfigure}[b]{0.13\textwidth} 
        \centering 
        \includegraphics[scale=0.16]{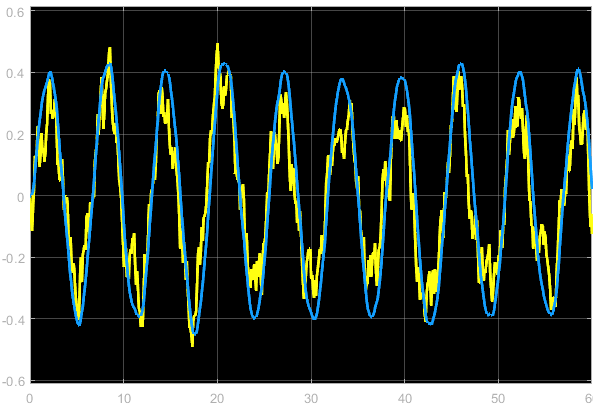} 
    \end{subfigure} 
    \begin{subfigure}[b]{0.13\textwidth}
        \centering 
        \includegraphics[scale=0.16]{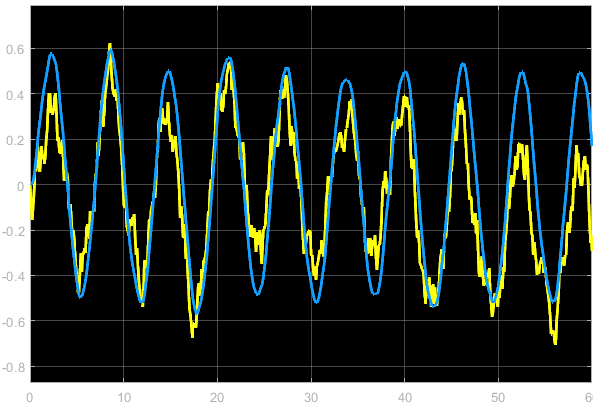} 
    \end{subfigure}
    \begin{subfigure}[b]{0.13\textwidth}
        \centering 
        \includegraphics[scale=0.16]{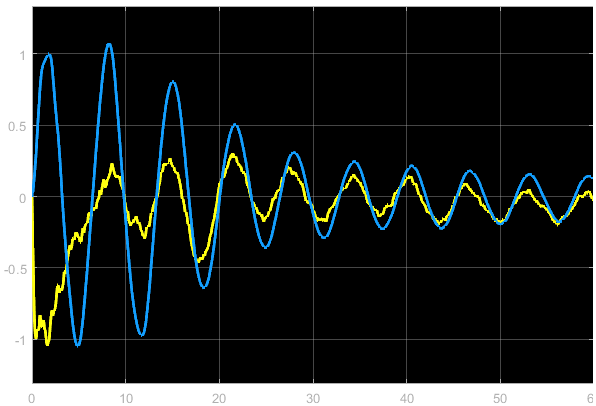} 
    \end{subfigure}
    \begin{subfigure}[b]{0.13\textwidth}
        \centering 
        \includegraphics[scale=0.16]{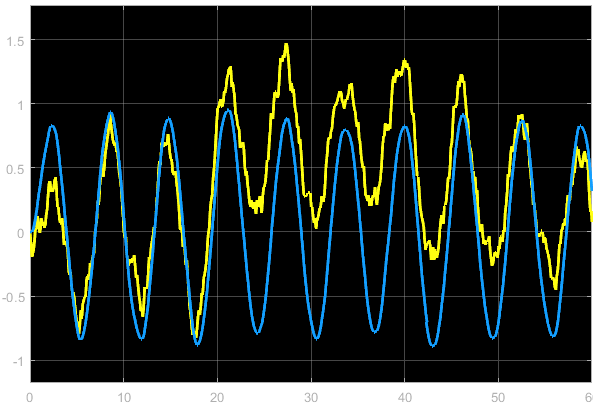} 
    \end{subfigure}
    \begin{subfigure}[b]{0.13\textwidth} 
        \centering 
        \includegraphics[scale=0.16]{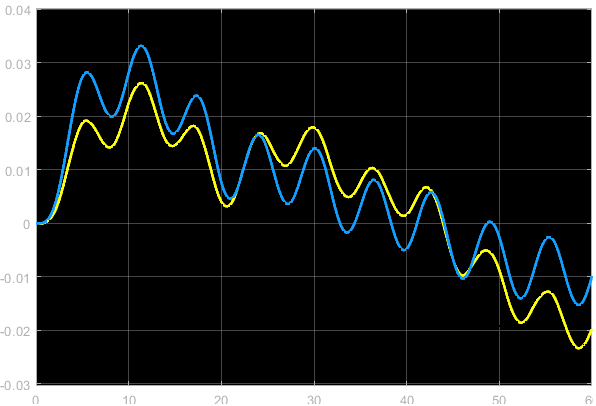} 
    \end{subfigure} 
    \begin{subfigure}[b]{0.13\textwidth} 
        \centering 
        \includegraphics[scale=0.16]{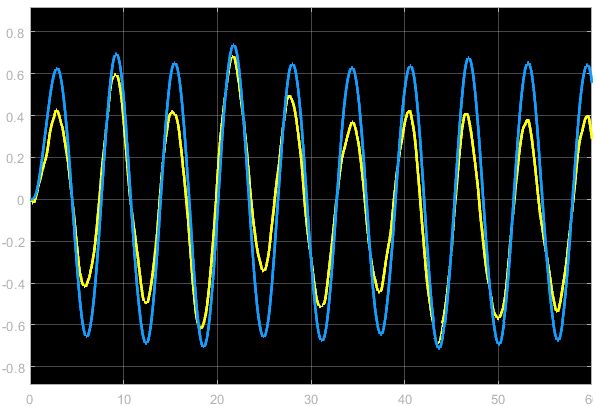} 
    \end{subfigure} 
    \begin{subfigure}[b]{0.13\textwidth} 
        \centering 
        \includegraphics[scale=0.16]{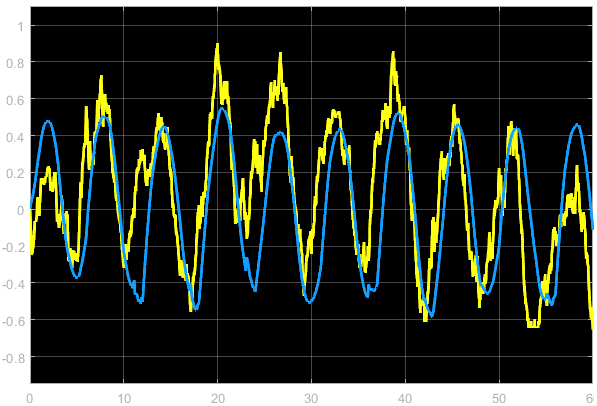} 
    \end{subfigure} 
    \caption{The system response of Classical GAN on various optimizers. Optimizer from left to right is respectively SGD, SGDM, Adam, PID, LPFSGD, HPFSGD, and FuzzyPID. The input of generator is a random Gaussian noise, and the input of discriminator is a square wave. (Blue wave comes from the discriminator, and yellow wave is the response of the generator.)} 
    \label{Figure2} 
\end{figure*}

GAN is designed to generate samples from the Gaussian noise. The performance of the GAN depends on its architecture \cite{zhou2023gan}. The generative network uses random inputs to generate samples, and the discriminative network aims to classify whether the generated sample can be classified \cite{2014Generative}. We present the control system of GAN in Figure \ref{Figure1} and its simulation in Figure \ref{Figure2}. We get its $\theta(s)$ as below:
\begin{equation} \label{Eq16}
\theta_D(s) = Opt \cdot \theta_{G}(s) \cdot E(s),
\end{equation}
\begin{equation} \label{Eq17}
\theta_G(s) = Opt \cdot E(s),
\end{equation} 
\begin{equation} \label{Eq18}
E(s) = \frac{\theta_{D}^{\ast}}{s} - \theta_{D}(s),
\end{equation} 
\noindent where $Opt$ is the optimization method, $\theta_D(s)$ is the desired Discriminator, $\theta_G(s)$ is the desired Generator. $E(s)$ is the feed-back error. $\theta_{G}^{\ast}$ is the optimal solution of the generator, and $\theta_{D}^{\ast}$ is the optimal solution of the discriminator. 

Eventually, we simplify $\theta_G(s)$ and $\theta_D(s)$ as below:
\begin{equation} \label{Eq19}
\theta_G(s) = \frac{1}{2} \cdot \left( \frac{\theta_{D}^{\ast}}{Opt} \pm \sqrt{(\frac{\theta_D^\ast}{Opt})^{2} - \frac{4}{s}} \right),
\end{equation}
\begin{equation} \label{Eq20}
\theta_D(s) = \theta^{2}_G(s),
\end{equation}
\noindent where if set $\theta_G(s)=0$, we get one pole point $s=0$. When using SGD as the $Opt$, $\theta_G(s)$ is a marginally stable system.

\section{The Laplace Transform of CycleGAN}
\begin{figure}[ht]
\centering
\includegraphics[angle=90,scale=0.46]{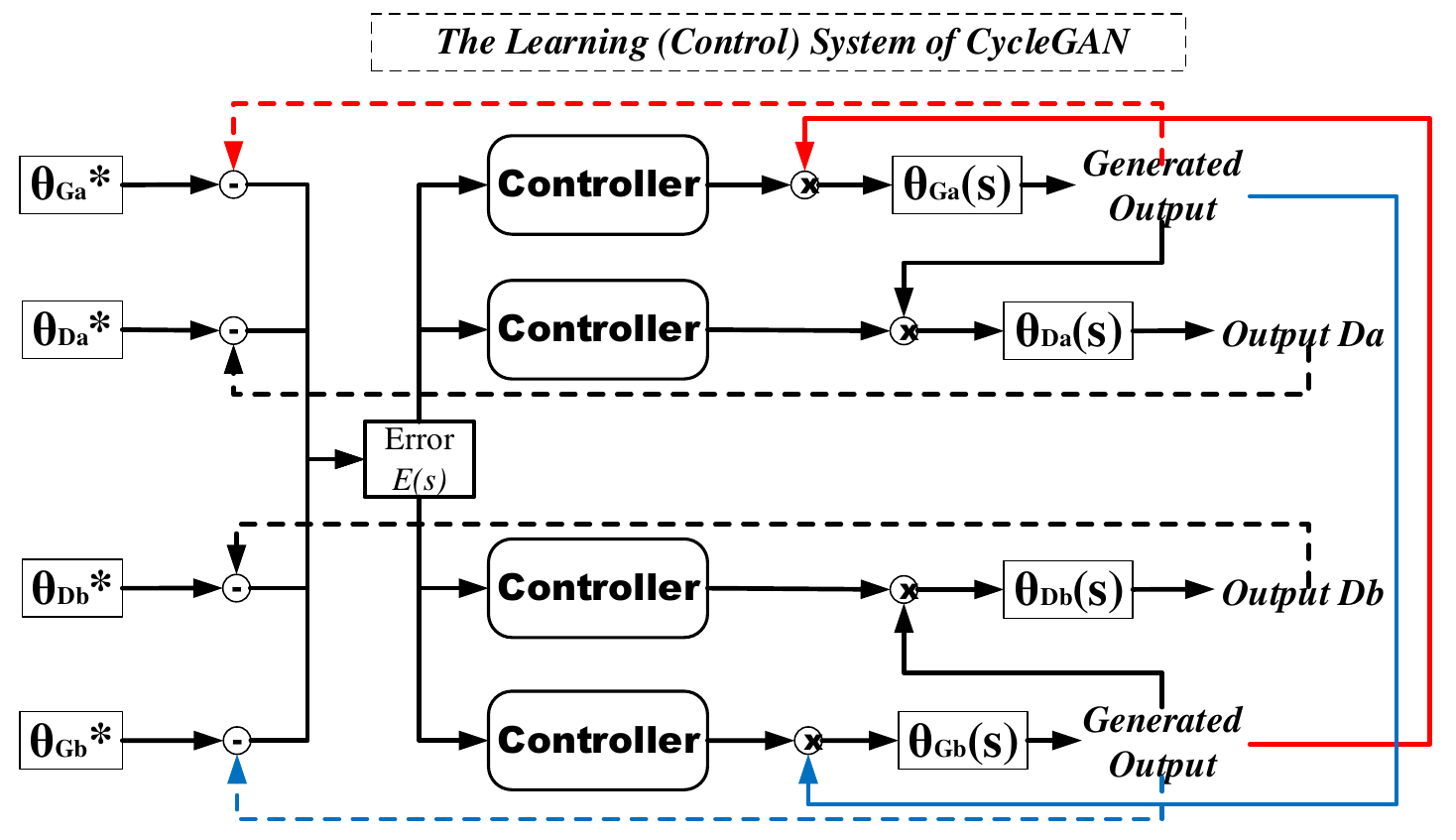}
\captionof{figure}{The control system of CycleGAN.}
\label{Figure3}
\end{figure}

\begin{figure*}[htp!]
    \centering 
    \begin{subfigure}[b]{0.24\textwidth} 
        \centering 
        \includegraphics[scale=0.28]{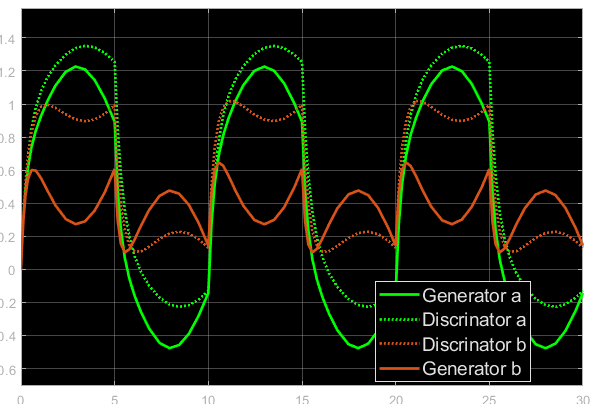} 
        \caption{\small CycleGAN on SGD.} 
        \label{Figure12-a} 
    \end{subfigure} 
    \begin{subfigure}[b]{0.24\textwidth}
        \centering 
        \includegraphics[scale=0.28]{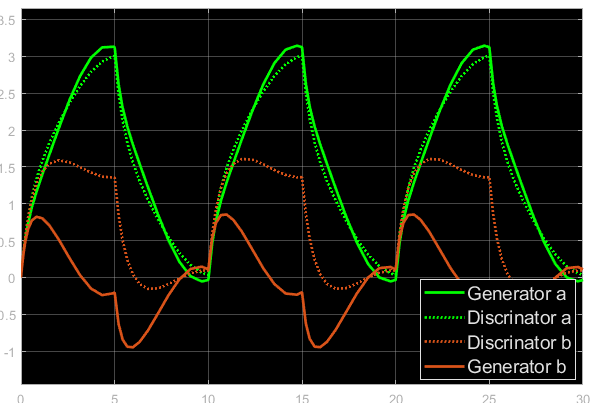} 
        \caption{\small CycleGAN on SGDM.} 
        \label{Figure12-b} 
    \end{subfigure}
    \begin{subfigure}[b]{0.24\textwidth}
        \centering 
        \includegraphics[scale=0.28]{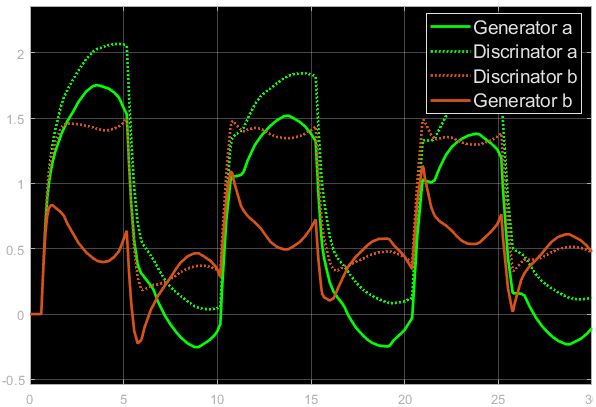} 
        \caption{\small  CycleGAN on Adam.} 
        \label{Figure12-c} 
    \end{subfigure}
    \begin{subfigure}[b]{0.24\textwidth}
        \centering 
        \includegraphics[scale=0.28]{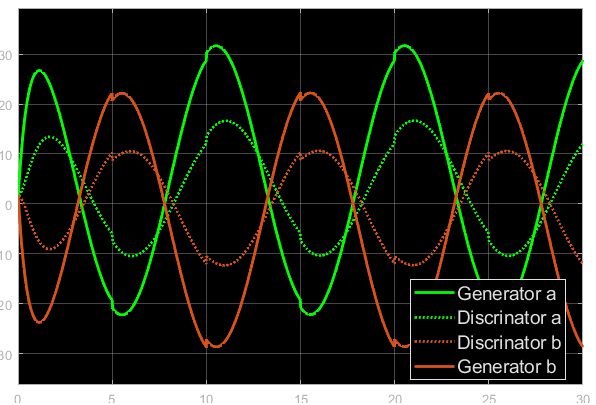} 
        \caption{\small  CycleGAN on PID.} 
        \label{Figure12-d} 
    \end{subfigure}
    \begin{subfigure}[b]{0.24\textwidth}
        \centering 
        \includegraphics[scale=0.28]{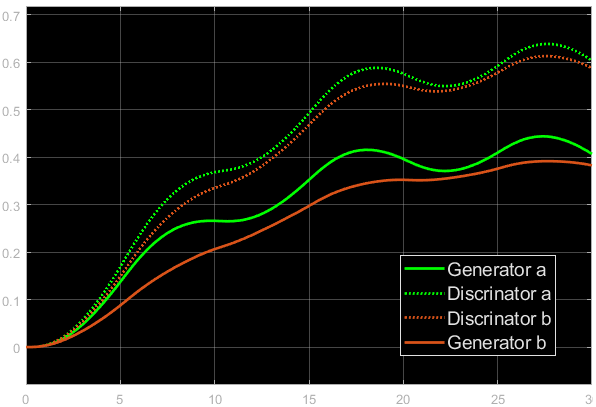} 
        \caption{\small  CycleGAN on LPF-SGD.} 
        \label{Figure12-e} 
    \end{subfigure}
    \begin{subfigure}[b]{0.24\textwidth}
        \centering 
        \includegraphics[scale=0.28]{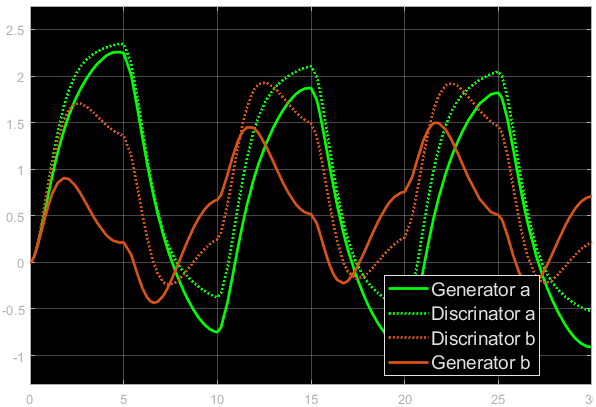} 
        \caption{\small  CycleGAN on HPF-SGD.} 
        \label{Figure12-f} 
    \end{subfigure}
        \begin{subfigure}[b]{0.24\textwidth}
        \centering 
        \includegraphics[scale=0.28]{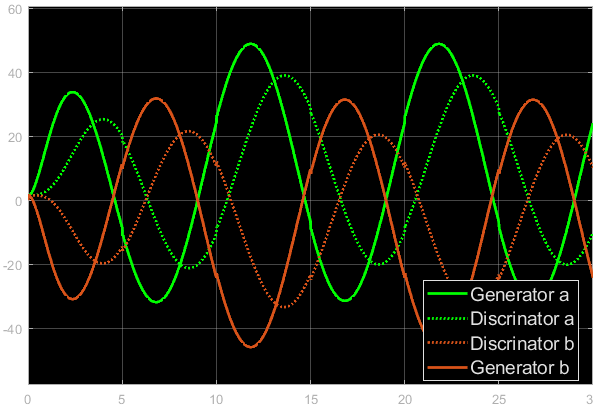} 
        \caption{\small CycleGAN on FuzzyPID.} 
        \label{Figure12-g} 
    \end{subfigure}
    \caption{The system response of CycleGAN on different optimizers, such as SGD, SGDM, Adam, PID, LPF-SGD, HPF-SGD and FuzzyPID optimizers. For $G_{a}$, the source is a sinusoidal wave (green), and the dashed green wave comes from $D_{a}$.  For $G_{a}$, the source is a cosine wave (red), and the dashed red wave comes from  $D_{b}$.} 
    \label{Figure4} 
\end{figure*}

CycleGAN \cite{zhu2017unpaired} aims to translate an image from a source domain $a$ to a target domain $b$ in the absence of paired examples. We denote the data distribution as $a \thicksim p_{data}(a)$ and $b \thicksim p_{data}(b)$. CycleGAN contains two mapping functions $G_{a}$: $A \rightarrow B$ and $G_{b}$: $B \rightarrow A$, and associated adversarial discriminators $D_{a}$ and $D_{b}$. {$D_{b}$ encourages generator $G_{a}$ to translate $A$ into outputs indistinguishable from domain $B$, and vice versa for $D_{A}$ and $B$. According to its learning system, we present the control system of CycleGAN as depicted in Figure \ref{Figure3} and its simulation as depicted in Figure \ref{Figure4}. CycleGAN has two generators and two discriminators. Meanwhile, ti takes the cycle consistency loss into account. The loss function has three parts as below: 
\begin{equation} \label{Eq53}
\begin{aligned}
\mathcal{L}(G_{a},G_{b},&D_{a},D_{b}) = \mathcal{L}_{GAN}(G_{a},D_{b},A,B) +\\
& \mathcal{L}_{GAN}(G_{b},D_{a},B,A) + \lambda \mathcal{L}_{cyc}(G_{a},G_{b}),
\end{aligned}
\end{equation}
where for the mapping function $G_{a}$: $A \rightarrow B$ and its discriminator $D_{b}$, we express the objective as:
\begin{equation} \label{Eq54}
\begin{aligned}
 \mathcal{L}_{GAN}(G_{a},D_{b},A,B) =& \mathbb{E}_{a \thicksim p_{data}(a)} [log D_{b} (b)] +\\
& \mathbb{E}_{b \thicksim p_{data}(b)}[log (1+D_{b}(G(a))]
\end{aligned}
\end{equation}

For each image $b$ from domain $B$ , $G_{a}$ and $G_{b}$ should satisfy the consistency of the backward cycle: $b \rightarrow G_{a}(b) \rightarrow G_{a}(G_{b}(b)) \approx y$. Thus, the cycle consistency loss is:

\begin{equation} \label{Eq55}
\begin{aligned}
 \mathcal{L}_{cyc}(G_{a},G_{b}) =& \mathbb{E}_{a \thicksim p_{data}(a)} [|| G_{b}(G_{a}(a))-a ||_{1}] +\\
& \mathbb{E}_{b \thicksim p_{data}(b)}[|| G_{a}(G_{b}(b))-b ||_{1}]
\end{aligned}
\end{equation}

CycleGAN used the L1 norm in this loss with an adversarial loss between $G_{b}(G_{a}(a))$ and $a$, and between $G_{a}(G_{b}(b))$ and $b$, but it did not observe the Laplace transform. Therefore, we get the system function of CycleGAN as below:

\begin{equation} \label{Eq56}
\begin{aligned}
\theta_{Da}(s) &= Opt \cdot \theta_{Ga}(s) \cdot E(s),
\\\theta_{Ga}(s) &= Opt \cdot E(s)    ,
\\\theta_{Db}(s) &= Opt \cdot \theta_{Gb}(s) \cdot E(s),
\\\theta_{Gb}(s) &= Opt \cdot E(s).
\end{aligned}
\end{equation} 
and
\begin{equation} \label{Eq60}
\begin{aligned}
E(s) = &\left[ \frac{\theta_{Da}^{\ast}}{s} - \theta_{Da}(s) \right] + \left[ \frac{\theta_{Db}^{\ast}}{s} - \theta_{Db}(s) \right] +\\
& \left[ \frac{\theta_{Ga}^{\ast}}{s} - \theta_{Ga}(s)\theta_{Gb}(s) \right] + \left[ \frac{\theta_{Gb}^{\ast}}{s} - \theta_{Gb}(s)\theta_{Ga}(s) \right].
\end{aligned}
\end{equation}

\section{Experiments}

\subsection{Experiment Settings} 
To analyse each optimizer, we formalise their system functions using the Laplace transform and simulate their system responses for analysing their attributes. We experiment these seven optimizers on MNIST \cite{lecun1998gradient} and UPSP \cite{carlucci2019hallucinating}.\footnote{Data is available at : \url{https://paperswithcode.com/sota/domain-adaptation-on-mnist-to-usps}.} As we believe that the training process of most GMs can be modeled as the source response of their control systems, we use Simulink (MATLAB R2022a) to simulate their response under different sources. All the experiments are implemented in PyTorch \cite{paszke2019pytorch} and conducted on a single NVIDIA A100 40GB GPU. All hyperparameters are presented in Table \ref{Table1} and Table \ref{Table2}.

To illustrate the connection between the control theory and the learning process of some complex artificial neural networks (ANNs), we analyze classical GAN \cite{2014Generative}, CycleGAN \cite{chu2017cyclegan}, and DDPM \cite{ho2020denoising} using Laplace transform. To verify the influence of different optimizers on GMs, we employ SGD, SGDM, Adam, PID, LPF-SGD, HPF-SGD and fuzzyPID on these three GMs to generate handwriting numbers.

\begin{table}[h]
\caption{Hyper-parameters for the image generation task on MNIST and UPSP.}
\label{Table1}
\scriptsize
\centering
\begin{tabular}{cccc}
\toprule
\textbf{Hyper-parameter}  & \textbf{\begin{tabular}[c]{@{}c@{}}GAN\end{tabular}} & \textbf{\begin{tabular}[c]{@{}c@{}}CycleGAN\end{tabular}} & \textbf{\begin{tabular}[c]{@{}c@{}}DDPM\end{tabular}}\\ \hline
Data Name                 & MNIST  & MNIST and UPSP & MNIST\\
Data augmentation                 & Auto  & Auto & Auto\\
Input resolution                  & [28,28,1]  & [28,28,1]  & [28,28,1]\\
Epochs                           & 200  & 200 & 100\\
Batch size                       & 200  & 200  & 200\\
Hidden dropout                    & 0  & 0 & 0\\
Random erasing prob               & 0 & 0 & 0\\
Random erasing prob               & 0 & 0 & 0\\
Random erasing prob               & 0  & 0 & 0\\
EMA decay                         & 0  & 0 & 0\\
Cutmix $\alpha$                   & 0  & 0 & 0\\
Mixup $\alpha$                    & 0  & 0 & 0\\
Cutmix-Mixup                      & 0 & 0 & 0\\
Peak learning rate                & 2e-4  & 2e-4  & 2e-4\\
\bottomrule
\end{tabular}
\end{table}

\begin{table}[h]
\caption{Coefficients of LPF-SGD and HPF-SGD using second-order IIR structure.}
\label{Table2}
\scriptsize
\centering
\begin{tabular}{cccccccc}
\toprule
\multirow{2}{*}{\textbf{Filter Type}}   & Gain  & \multicolumn{3}{c}{Numerator} & \multicolumn{3}{c}{Denominator} \\ \cline{2-8}
 & $G$ & $x_{0}$ & $x_{1}$  & $x_{2}$  & $y_{0}$  & $y_{1}$  & $y_{2}$\\ \hline
Low Pass Filter             & 0.49968  & 1 & -0.99937 & 0.00063 & 1.0  & 0  & -1.0   \\ 
High Pass Filter             & 0.49968  & 1 & 0.99937 & 0.00063 & 1.0  & 0  & -1.0   \\
\bottomrule
\end{tabular}
\end{table}

\subsection{Simulation} 

To get a clear analysis result, we merely use a random noise source and a sinusoidal wave source for the generator. Additionally, we apply a square wave source to the discriminator. As we are not able to know the real system function of each GM, we roughly set their system order as 2.

\section{Results and Analysis}

\subsection{Simulation Performance}

\begin{figure}
    \centering 
    \begin{subfigure}[b]{0.24\textwidth} 
        \caption{\tiny $epoch_{1}$.} 
        \centering 
        \includegraphics[scale=0.18]{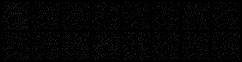} 
    \end{subfigure} 
    \begin{subfigure}[b]{0.24\textwidth} 
        \caption{\tiny $epoch_{50}$.} 
        \centering 
        \includegraphics[scale=0.18]{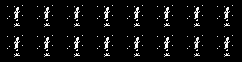}
    \end{subfigure}
    \begin{subfigure}[b]{0.24\textwidth} 
        \caption{\tiny $epoch_{100}$.}
        \centering 
        \includegraphics[scale=0.18]{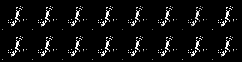}
    \end{subfigure} 
    \begin{subfigure}[b]{0.24\textwidth} 
        \caption{\tiny $epoch_{200}$.}
        \centering 
        \includegraphics[scale=0.18]{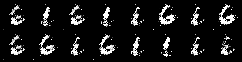}
    \end{subfigure} 
    \vskip 
    \baselineskip 
    \begin{subfigure}[b]{0.24\textwidth} 
        \centering 
        \includegraphics[scale=0.18]{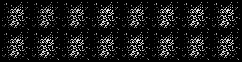} 
    \end{subfigure} 
    \begin{subfigure}[b]{0.24\textwidth} 
        \centering 
        \includegraphics[scale=0.18]{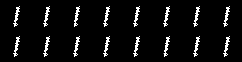} 
    \end{subfigure}
    \begin{subfigure}[b]{0.24\textwidth} 
        \centering 
        \includegraphics[scale=0.18]{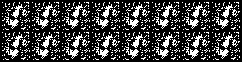} 
    \end{subfigure} 
    \begin{subfigure}[b]{0.24\textwidth} 
        \centering 
        \includegraphics[scale=0.18]{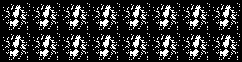} 
    \end{subfigure} 
    \vskip 
    \baselineskip 
    \begin{subfigure}[b]{0.24\textwidth} 
        \centering 
        \includegraphics[scale=0.18]{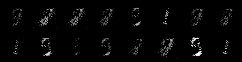} 
    \end{subfigure} 
    \begin{subfigure}[b]{0.24\textwidth} 
        \centering 
        \includegraphics[scale=0.18]{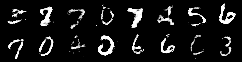} 
    \end{subfigure}
    \begin{subfigure}[b]{0.24\textwidth} 
        \centering 
        \includegraphics[scale=0.18]{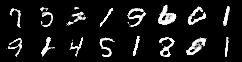} 
    \end{subfigure} 
    \begin{subfigure}[b]{0.24\textwidth} 
        \centering 
        \includegraphics[scale=0.18]{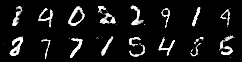} 
    \end{subfigure} 
    \vskip 
    \baselineskip 
    \begin{subfigure}[b]{0.24\textwidth} 
        \centering 
        \includegraphics[scale=0.18]{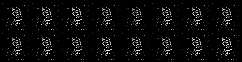} 
    \end{subfigure} 
    \begin{subfigure}[b]{0.24\textwidth} 
        \centering 
        \includegraphics[scale=0.18]{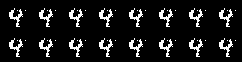} 
    \end{subfigure}
    \begin{subfigure}[b]{0.24\textwidth} 
        \centering 
        \includegraphics[scale=0.18]{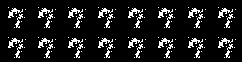} 
    \end{subfigure} 
    \begin{subfigure}[b]{0.24\textwidth} 
        \centering 
        \includegraphics[scale=0.18]{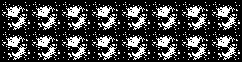} 
    \end{subfigure}\
    \vskip 
    \baselineskip 
    \begin{subfigure}[b]{0.24\textwidth} 
        \centering 
        \includegraphics[scale=0.18]{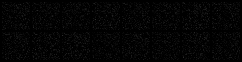} 
    \end{subfigure} 
    \begin{subfigure}[b]{0.24\textwidth} 
        \centering 
        \includegraphics[scale=0.18]{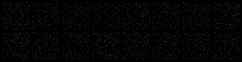} 
    \end{subfigure}
    \begin{subfigure}[b]{0.24\textwidth} 
        \centering 
        \includegraphics[scale=0.18]{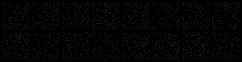} 
    \end{subfigure} 
    \begin{subfigure}[b]{0.24\textwidth} 
        \centering 
        \includegraphics[scale=0.18]{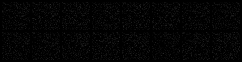} 
    \end{subfigure} 
    \vskip
    \baselineskip 
    \begin{subfigure}[b]{0.24\textwidth} 
        \centering 
        \includegraphics[scale=0.18]{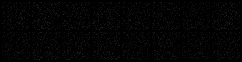} 
    \end{subfigure} 
    \begin{subfigure}[b]{0.24\textwidth} 
        \centering 
        \includegraphics[scale=0.18]{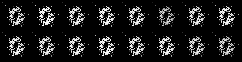} 
    \end{subfigure}
    \begin{subfigure}[b]{0.24\textwidth} 
        \centering 
        \includegraphics[scale=0.18]{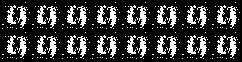} 
    \end{subfigure} 
    \begin{subfigure}[b]{0.24\textwidth} 
        \centering 
        \includegraphics[scale=0.18]{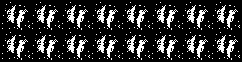} 
    \end{subfigure} 
    \vskip 
    \baselineskip 
    \begin{subfigure}[b]{0.24\textwidth} 
        \centering 
        \includegraphics[scale=0.18]{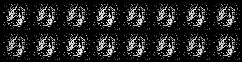} 
    \end{subfigure} 
    \begin{subfigure}[b]{0.24\textwidth} 
        \centering 
        \includegraphics[scale=0.18]{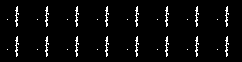} 
    \end{subfigure}
    \begin{subfigure}[b]{0.24\textwidth} 
        \centering 
        \includegraphics[scale=0.18]{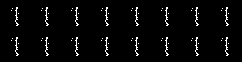} 
    \end{subfigure} 
    \begin{subfigure}[b]{0.24\textwidth} 
        \centering 
        \includegraphics[scale=0.18]{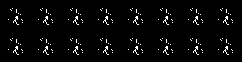} 
    \end{subfigure} 
    \caption{The generated samples from classical GAN on corresponding optimizers (from top to bottom is respectively SGD, SGDM, Adam, PID, LPF-SGD, HPF-SGD, and FuzzyPID).} 
    \label{Figure5} 
\end{figure}

\begin{figure*}[htp!]
    \centering 
    \begin{subfigure}[b]{0.1\textwidth}
        \caption{\scriptsize $G_{a}$ to $G_{b}$ on $epoch_{1}$.} 
        \centering 
        \includegraphics[scale=0.40]{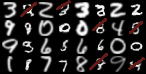} 
    \end{subfigure} 
    \hfill 
    \begin{subfigure}[b]{0.1\textwidth}
        \caption{\scriptsize $G_{a}$ to $G_{b}$ on $epoch_{50}$.} 
        \centering 
        \includegraphics[scale=0.40]{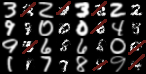} 
    \end{subfigure}
    \hfill
    \begin{subfigure}[b]{0.1\textwidth} 
        \caption{\scriptsize $G_{a}$ to $G_{b}$ on $epoch_{100}$.} 
        \centering 
        \includegraphics[scale=0.40]{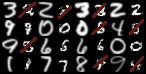} 
    \end{subfigure} 
    \hfill 
    \begin{subfigure}[b]{0.1\textwidth} 
        \caption{\scriptsize $G_{a}$ to $G_{b}$ on $epoch_{200}$.} 
        \centering 
        \includegraphics[scale=0.40]{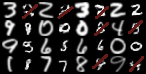} 
    \end{subfigure}
    \hfill 
    \begin{subfigure}[b]{0.1\textwidth} 
        \caption{\scriptsize $G_{b}$ to $G_{a}$ on $epoch_{1}$.} 
        \centering 
        \includegraphics[scale=0.40]{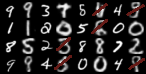} 
    \end{subfigure} 
    \hfill 
    \begin{subfigure}[b]{0.1\textwidth} 
        \caption{\scriptsize $G_{b}$ to $G_{a}$ on $epoch_{50}$.} 
        \centering 
        \includegraphics[scale=0.40]{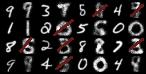} 
    \end{subfigure}
    \hfill
    \begin{subfigure}[b]{0.1\textwidth} 
        \caption{\scriptsize $G_{b}$ to $G_{a}$ on $epoch_{100}$.} 
        \centering 
        \includegraphics[scale=0.40]{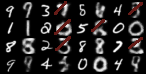} 
    \end{subfigure} 
    \hfill 
    \begin{subfigure}[b]{0.1\textwidth} 
        \caption{\scriptsize $G_{b}$ to $G_{a}$ on $epoch_{200}$.} 
        \centering 
        \includegraphics[scale=0.40]{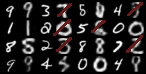} 
    \end{subfigure} 
    \vskip 
    \baselineskip  
    \begin{subfigure}[b]{0.1\textwidth} 
        \centering 
        \includegraphics[scale=0.40]{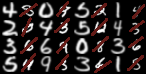} 
    \end{subfigure} 
    \hfill 
    \begin{subfigure}[b]{0.1\textwidth} 
        \centering 
        \includegraphics[scale=0.40]{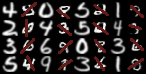} 
    \end{subfigure}
    \hfill
    \begin{subfigure}[b]{0.1\textwidth} 
        \centering 
        \includegraphics[scale=0.40]{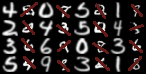} 
    \end{subfigure} 
    \hfill 
    \begin{subfigure}[b]{0.1\textwidth} 
        \centering 
        \includegraphics[scale=0.40]{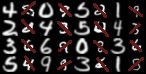} 
    \end{subfigure}
    \hfill 
    \begin{subfigure}[b]{0.1\textwidth} 
        \centering 
        \includegraphics[scale=0.40]{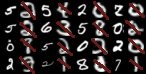} 
    \end{subfigure} 
    \hfill 
    \begin{subfigure}[b]{0.1\textwidth} 
        \centering 
        \includegraphics[scale=0.40]{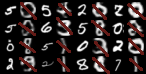}
    \end{subfigure}
    \hfill
    \begin{subfigure}[b]{0.1\textwidth} 
        \centering 
        \includegraphics[scale=0.40]{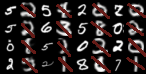} 
    \end{subfigure} 
    \hfill 
    \begin{subfigure}[b]{0.1\textwidth} 
        \centering 
        \includegraphics[scale=0.40]{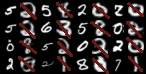} 
    \end{subfigure} 
    \vskip 
    \baselineskip  
    \begin{subfigure}[b]{0.1\textwidth} 
        \centering 
        \includegraphics[scale=0.40]{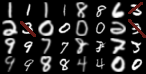} 
    \end{subfigure} 
    \hfill 
    \begin{subfigure}[b]{0.1\textwidth} 
        \centering 
        \includegraphics[scale=0.40]{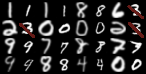} 
    \end{subfigure}
    \hfill
    \begin{subfigure}[b]{0.1\textwidth} 
        \centering 
        \includegraphics[scale=0.40]{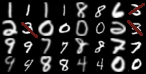} 
    \end{subfigure} 
    \hfill 
    \begin{subfigure}[b]{0.1\textwidth} 
        \centering 
        \includegraphics[scale=0.40]{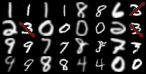} 
    \end{subfigure}
    \hfill 
    \begin{subfigure}[b]{0.1\textwidth} 
        \centering 
        \includegraphics[scale=0.40]{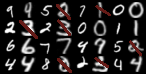} 
    \end{subfigure} 
    \hfill 
    \begin{subfigure}[b]{0.1\textwidth} 
        \centering 
        \includegraphics[scale=0.40]{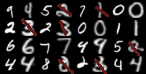} 
    \end{subfigure}
    \hfill
    \begin{subfigure}[b]{0.1\textwidth} 
        \centering 
        \includegraphics[scale=0.40]{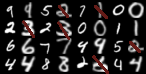} 
    \end{subfigure} 
    \hfill 
    \begin{subfigure}[b]{0.1\textwidth} 
        \centering 
        \includegraphics[scale=0.40]{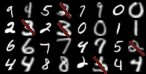} 
    \end{subfigure}
    \vskip 
    \baselineskip  
    \begin{subfigure}[b]{0.1\textwidth} 
        \centering 
        \includegraphics[scale=0.40]{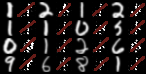} 
    \end{subfigure} 
    \hfill 
    \begin{subfigure}[b]{0.1\textwidth} 
        \centering 
        \includegraphics[scale=0.40]{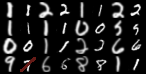} 
    \end{subfigure}
    \hfill
    \begin{subfigure}[b]{0.1\textwidth} 
        \centering 
        \includegraphics[scale=0.40]{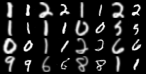} 
    \end{subfigure} 
    \hfill 
    \begin{subfigure}[b]{0.1\textwidth} 
        \centering 
        \includegraphics[scale=0.40]{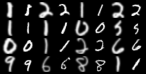} 
    \end{subfigure}
    \hfill 
    \begin{subfigure}[b]{0.1\textwidth} 
        \centering 
        \includegraphics[scale=0.40]{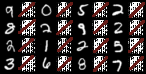} 
    \end{subfigure} 
    \hfill 
    \begin{subfigure}[b]{0.1\textwidth} 
        \centering 
        \includegraphics[scale=0.40]{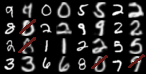} 
    \end{subfigure}
    \hfill
    \begin{subfigure}[b]{0.1\textwidth} 
        \centering 
        \includegraphics[scale=0.40]{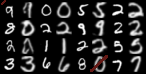} 
    \end{subfigure} 
    \hfill 
    \begin{subfigure}[b]{0.1\textwidth} 
        \centering 
        \includegraphics[scale=0.40]{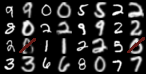} 
    \end{subfigure}
            \vskip 
    \baselineskip  
    \begin{subfigure}[b]{0.1\textwidth} 
        \centering 
        \includegraphics[scale=0.40]{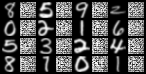} 
    \end{subfigure} 
    \hfill 
    \begin{subfigure}[b]{0.1\textwidth} 
        \centering 
        \includegraphics[scale=0.40]{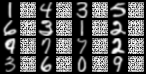} 
    \end{subfigure}
    \hfill
    \begin{subfigure}[b]{0.1\textwidth} 
        \centering 
        \includegraphics[scale=0.40]{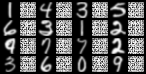} 
    \end{subfigure} 
    \hfill 
    \begin{subfigure}[b]{0.1\textwidth} 
        \centering 
        \includegraphics[scale=0.40]{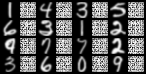} 
    \end{subfigure}
    \hfill 
    \begin{subfigure}[b]{0.1\textwidth} 
        \centering 
        \includegraphics[scale=0.40]{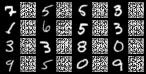} 
    \end{subfigure} 
    \hfill 
    \begin{subfigure}[b]{0.1\textwidth} 
        \centering 
        \includegraphics[scale=0.40]{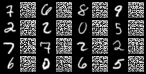} 
    \end{subfigure}
    \hfill
    \begin{subfigure}[b]{0.1\textwidth} 
        \centering 
        \includegraphics[scale=0.40]{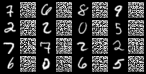} 
    \end{subfigure} 
    \hfill 
    \begin{subfigure}[b]{0.1\textwidth} 
        \centering 
        \includegraphics[scale=0.40]{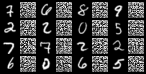} 
    \end{subfigure}
    \vskip
    \baselineskip  
    \begin{subfigure}[b]{0.1\textwidth} 
        \centering 
        \includegraphics[scale=0.40]{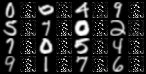} 
    \end{subfigure} 
    \hfill 
    \begin{subfigure}[b]{0.1\textwidth} 
        \centering 
        \includegraphics[scale=0.40]{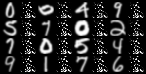} 
    \end{subfigure}
    \hfill
    \begin{subfigure}[b]{0.1\textwidth} 
        \centering 
        \includegraphics[scale=0.40]{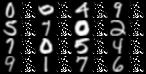} 
    \end{subfigure} 
    \hfill 
    \begin{subfigure}[b]{0.1\textwidth} 
        \centering 
        \includegraphics[scale=0.40]{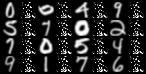} 
    \end{subfigure}
    \hfill 
    \begin{subfigure}[b]{0.1\textwidth} 
        \centering 
        \includegraphics[scale=0.40]{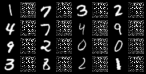} 
    \end{subfigure} 
    \hfill 
    \begin{subfigure}[b]{0.1\textwidth} 
        \centering 
        \includegraphics[scale=0.40]{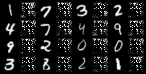} 
    \end{subfigure}
    \hfill
    \begin{subfigure}[b]{0.1\textwidth} 
        \centering 
        \includegraphics[scale=0.40]{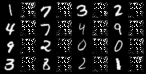} 
    \end{subfigure} 
    \hfill 
    \begin{subfigure}[b]{0.1\textwidth} 
        \centering 
        \includegraphics[scale=0.40]{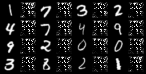} 
    \end{subfigure}
    \vskip 
    \baselineskip  
    \begin{subfigure}[b]{0.1\textwidth} 
        \centering 
        \includegraphics[scale=0.40]{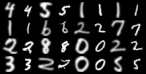} 
    \end{subfigure} 
    \hfill 
    \begin{subfigure}[b]{0.1\textwidth} 
        \centering 
        \includegraphics[scale=0.40]{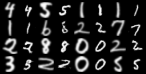} 
    \end{subfigure}
    \hfill
    \begin{subfigure}[b]{0.1\textwidth} 
        \centering 
        \includegraphics[scale=0.40]{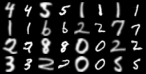} 
    \end{subfigure} 
    \hfill 
    \begin{subfigure}[b]{0.1\textwidth} 
        \centering 
        \includegraphics[scale=0.40]{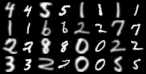} 
    \end{subfigure}
    \hfill 
    \begin{subfigure}[b]{0.1\textwidth} 
        \centering 
        \includegraphics[scale=0.40]{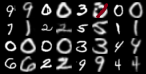} 
    \end{subfigure} 
    \hfill 
    \begin{subfigure}[b]{0.1\textwidth} 
        \centering 
        \includegraphics[scale=0.40]{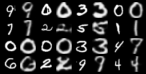} 
    \end{subfigure}
    \hfill
    \begin{subfigure}[b]{0.1\textwidth} 
        \centering 
        \includegraphics[scale=0.40]{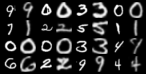} 
    \end{subfigure} 
    \hfill 
    \begin{subfigure}[b]{0.1\textwidth} 
        \centering 
        \includegraphics[scale=0.40]{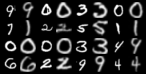} 
    \end{subfigure}
    \caption{The generated samples from CycleGAN on corresponding optimizers (From top to bottom is SGD, SGDM, Adam, PID, LPF-SGD, HPF-SGD and FuzzyPID). The generated confident errors are marked with a red line.} 
    \label{Figure6} 
\end{figure*}

We simulate the system response of GAN and CycleGAN using Simulink over seven optimizers. Adam outperforms other optimizers when using GAN to generate samples from random noise (seen in Figure \ref{Figure2}). Adam is able to converge, and it gradually becomes stable. Although other optimizers (e.g. SGD, SGDM, LPF-SGD, PID, and FuzzyPID) do not have an obvious drift phenomenon at the early stage and meanwhile keep a marginal stability, at the late stage, they (SGDM, LPF-SGD, PID, and FuzzyPID) show an obvious baseline drift phenomenon.

When using CycleGAN to translate one image to another image, FuzzyPID optimizer presents the best simulation performance (seen in Figure \ref{Figure4}). PID optimizer also can generate a remarkable sinusoidal wave both on $G_{a}$ and $G_{b}$, but there are some breakages near each wave peak and bottom. Although SGD and AdaM can generate some discernible sinusoidal patterns, they have an increase of divergence risk. SGDM optimization methods fail to generate sinusoidal signals.

\subsection{Generated Samples of GAN and CycleGAN}

For the generated samples of GAN (seen in Figure \ref{Figure5}) and DDPM (seen in Figure \ref{Figure6}, \ref{Figure7}, \ref{Figure8}, and \ref{Figure9}), most optimizers can generate acceptable samples at the early generation stage, despite that apart from Adam, other optimizers only generate the duplicate sample. At the late generation stage, Adam still shows the dominant advantage, but other optimizers are not able to generate a recognisable sample, as their generated samples contain strong noise.

For the generated MNIST and UPSP in Figure \ref{Figure6}, after 100 epochs training, PID and FuzzPID can generate correct samples both from $G_{a}$ to $G_{b}$ and from $G_{b}$ to $G_{a}$. Notably, the ability of CycleGAN to produce samples from a single dataset was significantly enhanced when utilizing the FuzzyPID, which yielded flawless samples from the outset. This suggests that FuzzyPID might be the optimal choice for optimizing the learning updates of CycleGAN. A manual evaluation of the alignment between samples from and vice versa was also conducted. Preliminary observations indicate that the PID and FuzzyPID optimizers outshine the others when applied to models that utilize a cycle consistency loss, such as CycleGAN.

\subsection{Denoising Diffusion Probabilistic Models}

\begin{figure*}[h]
    \centering 
    \begin{subfigure}[b]{0.1\textwidth}
        \caption{\scriptsize SGD on $epoch_{1}$.} 
        \centering 
        \includegraphics[scale=0.34]{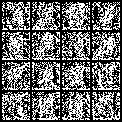} 
    \end{subfigure} 
    \hfill 
    \begin{subfigure}[b]{0.1\textwidth}
        \caption{\scriptsize SGD on $epoch_{25}$.} 
        \centering 
        \includegraphics[scale=0.34]{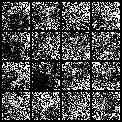} 
    \end{subfigure}
    \hfill
    \begin{subfigure}[b]{0.1\textwidth} 
        \caption{\scriptsize SGD on $epoch_{50}$.} 
        \centering 
        \includegraphics[scale=0.34]{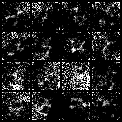} 
    \end{subfigure} 
    \hfill 
    \begin{subfigure}[b]{0.1\textwidth} 
        \caption{\scriptsize SGD on $epoch_{100}$.} 
        \centering 
        \includegraphics[scale=0.34]{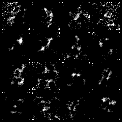} 
    \end{subfigure}
    \hfill 
    \begin{subfigure}[b]{0.1\textwidth} 
        \caption{\scriptsize SGDM on $epoch_{1}$.} 
        \centering 
        \includegraphics[scale=0.34]{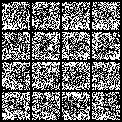} 
    \end{subfigure} 
    \hfill 
    \begin{subfigure}[b]{0.1\textwidth} 
        \caption{\scriptsize SGDM on $epoch_{25}$.} 
        \centering 
        \includegraphics[scale=0.34]{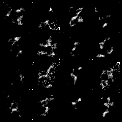} 
    \end{subfigure}
    \hfill
    \begin{subfigure}[b]{0.1\textwidth} 
        \caption{\scriptsize SGDM on $epoch_{50}$.} 
        \centering 
        \includegraphics[scale=0.34]{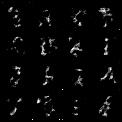} 
    \end{subfigure} 
    \hfill 
    \begin{subfigure}[b]{0.1\textwidth} 
        \caption{\scriptsize SGDM on $epoch_{100}$.} 
        \centering 
        \includegraphics[scale=0.34]{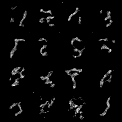} 
    \end{subfigure}
    \vskip 
    \baselineskip 
    \begin{subfigure}[b]{0.1\textwidth} 
        \caption{\scriptsize PID on $epoch_{1}$.} 
        \centering 
        \includegraphics[scale=0.34]{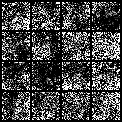} 
    \end{subfigure} 
    \hfill 
    \begin{subfigure}[b]{0.1\textwidth} 
        \caption{\scriptsize PID on $epoch_{25}$.} 
        \centering 
        \includegraphics[scale=0.34]{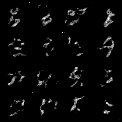} 
    \end{subfigure}
    \hfill
    \begin{subfigure}[b]{0.1\textwidth} 
        \caption{\scriptsize PID on $epoch_{50}$.} 
        \centering 
        \includegraphics[scale=0.34]{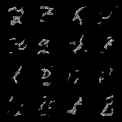} 
    \end{subfigure} 
    \hfill 
    \begin{subfigure}[b]{0.1\textwidth} 
        \caption{\scriptsize PID on $epoch_{100}$.} 
        \centering 
        \includegraphics[scale=0.34]{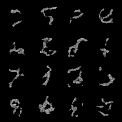} 
    \end{subfigure}
    \hfill 
    \begin{subfigure}[b]{0.1\textwidth} 
        \caption{\scriptsize FuzzyPID on $epoch_{1}$.} 
        \centering 
        \includegraphics[scale=0.34]{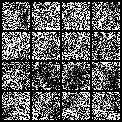} 
    \end{subfigure} 
    \hfill 
    \begin{subfigure}[b]{0.1\textwidth} 
        \caption{\scriptsize FuzzyPID on $epoch_{25}$.} 
        \centering 
        \includegraphics[scale=0.34]{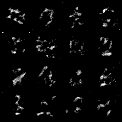} 
    \end{subfigure}
    \hfill
    \begin{subfigure}[b]{0.1\textwidth} 
        \caption{\scriptsize FuzzyPID on $epoch_{50}$.} 
        \centering 
        \includegraphics[scale=0.34]{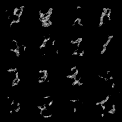} 
    \end{subfigure} 
    \hfill 
    \begin{subfigure}[b]{0.1\textwidth} 
        \caption{\scriptsize FuzzyPID on $epoch_{100}$.} 
        \centering 
        \includegraphics[scale=0.34]{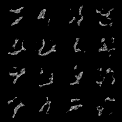} 
    \end{subfigure}
    \caption{The generated samples from DDPM on corresponding optimizers (From top to bottom is SGD, SGDM, Adam, PID, LPF-SGD, HPF-SGD and FuzzyPID).} 
    \label{Figure7} 
\end{figure*}

\begin{figure*}[h]
    \centering 
    \begin{subfigure}[b]{0.1\textwidth} 
        \caption{\scriptsize Adam on $epoch_{1}$.} 
        \centering 
        \includegraphics[scale=0.34]{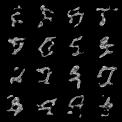} 
    \end{subfigure} 
    \hfill 
    \begin{subfigure}[b]{0.1\textwidth} 
        \caption{\scriptsize Adam on $epoch_{25}$.} 
        \centering 
        \includegraphics[scale=0.34]{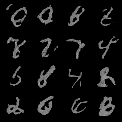} 
    \end{subfigure}
    \hfill
    \begin{subfigure}[b]{0.1\textwidth} 
        \caption{\scriptsize Adam on $epoch_{50}$.} 
        \centering 
        \includegraphics[scale=0.34]{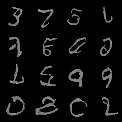} 
    \end{subfigure} 
    \hfill 
    \begin{subfigure}[b]{0.1\textwidth} 
        \caption{\scriptsize Adam on $epoch_{100}$.} 
        \centering 
        \includegraphics[scale=0.34]{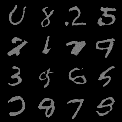} 
    \end{subfigure}
    \hfill 
    \begin{subfigure}[b]{0.1\textwidth} 
        \caption{\scriptsize AdamW on $epoch_{1}$.} 
        \centering 
        \includegraphics[scale=0.34]{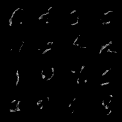} 
    \end{subfigure} 
    \hfill 
    \begin{subfigure}[b]{0.1\textwidth} 
        \caption{\scriptsize AdamW on $epoch_{25}$.} 
        \centering 
        \includegraphics[scale=0.34]{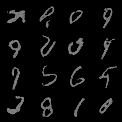} 
    \end{subfigure}
    \hfill
    \begin{subfigure}[b]{0.1\textwidth} 
        \caption{\scriptsize AdamW on $epoch_{50}$.} 
        \centering 
        \includegraphics[scale=0.34]{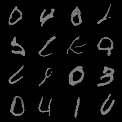} 
    \end{subfigure} 
    \hfill 
    \begin{subfigure}[b]{0.1\textwidth} 
        \caption{\scriptsize AdamW on $epoch_{100}$.} 
        \centering 
        \includegraphics[scale=0.34]{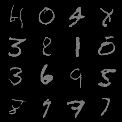} 
    \end{subfigure} 
    \vskip 
    \baselineskip  
    \begin{subfigure}[b]{0.1\textwidth} 
        \caption{\scriptsize RAdam on $epoch_{1}$.} 
        \centering 
        \includegraphics[scale=0.34]{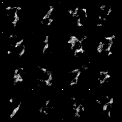} 
    \end{subfigure} 
    \hfill 
    \begin{subfigure}[b]{0.1\textwidth} 
        \caption{\scriptsize RAdam on $epoch_{25}$.} 
        \centering 
        \includegraphics[scale=0.34]{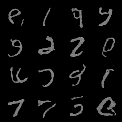} 
    \end{subfigure}
    \hfill
    \begin{subfigure}[b]{0.1\textwidth} 
        \caption{\scriptsize RAdam on $epoch_{50}$.} 
        \centering 
        \includegraphics[scale=0.34]{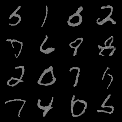} 
    \end{subfigure} 
    \hfill 
    \begin{subfigure}[b]{0.1\textwidth} 
        \caption{\scriptsize RAdam on $epoch_{100}$.} 
        \centering 
        \includegraphics[scale=0.34]{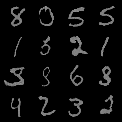} 
    \end{subfigure}
    \hfill 
    \begin{subfigure}[b]{0.1\textwidth} 
        \caption{\scriptsize NAdam on $epoch_{1}$.} 
        \centering 
        \includegraphics[scale=0.34]{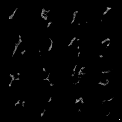} 
    \end{subfigure} 
    \hfill 
    \begin{subfigure}[b]{0.1\textwidth} 
        \caption{\scriptsize NAdam on $epoch_{25}$.} 
        \centering 
        \includegraphics[scale=0.34]{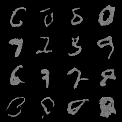} 
    \end{subfigure}
    \hfill
    \begin{subfigure}[b]{0.1\textwidth} 
        \caption{\scriptsize NAdam on $epoch_{50}$.} 
        \centering 
        \includegraphics[scale=0.34]{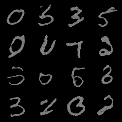} 
    \end{subfigure} 
    \hfill 
    \begin{subfigure}[b]{0.1\textwidth} 
        \caption{\scriptsize NAdam on $epoch_{100}$.} 
        \centering 
        \includegraphics[scale=0.34]{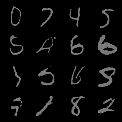} 
    \end{subfigure}

    \caption{The generated samples from DDPM on corresponding optimizers (From top to bottom is SGD, SGDM, Adam, PID, LPF-SGD, HPF-SGD and FuzzyPID).} 
    \label{Figure8} 
\end{figure*}

\begin{figure*}[h]
    \centering 
    \begin{subfigure}[b]{0.1\textwidth} 
        \caption{\scriptsize LPFSGD on $epoch_{1}$.} 
        \centering 
        \includegraphics[scale=0.34]{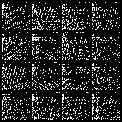} 
    \end{subfigure} 
    \hfill 
    \begin{subfigure}[b]{0.1\textwidth} 
        \caption{\scriptsize LPFSGD on $epoch_{25}$.} 
        \centering 
        \includegraphics[scale=0.34]{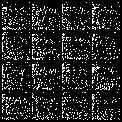} 
    \end{subfigure}
    \hfill
    \begin{subfigure}[b]{0.1\textwidth} 
        \caption{\scriptsize LPFSGD on $epoch_{50}$.} 
        \centering 
        \includegraphics[scale=0.34]{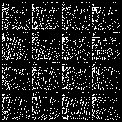} 
    \end{subfigure} 
    \hfill 
    \begin{subfigure}[b]{0.1\textwidth} 
        \caption{\scriptsize LPFSGD on $epoch_{100}$.} 
        \centering 
        \includegraphics[scale=0.34]{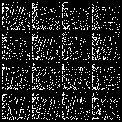} 
    \end{subfigure}
    \hfill 
    \begin{subfigure}[b]{0.1\textwidth} 
        \caption{\scriptsize HPFSGD on $epoch_{1}$.} 
        \centering 
        \includegraphics[scale=0.34]{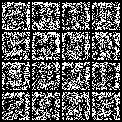} 
    \end{subfigure} 
    \hfill 
    \begin{subfigure}[b]{0.1\textwidth} 
        \caption{\scriptsize HPFSGD on $epoch_{25}$.} 
        \centering 
        \includegraphics[scale=0.34]{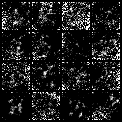} 
    \end{subfigure}
    \hfill
    \begin{subfigure}[b]{0.1\textwidth} 
        \caption{\scriptsize HPFSGD on $epoch_{50}$.} 
        \centering 
        \includegraphics[scale=0.34]{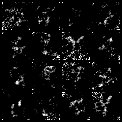} 
    \end{subfigure} 
    \hfill 
    \begin{subfigure}[b]{0.1\textwidth} 
        \caption{\scriptsize HPFSGD on $epoch_{100}$.} 
        \centering 
        \includegraphics[scale=0.34]{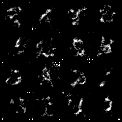} 
    \end{subfigure}
    \caption{The generated samples from DDPM on corresponding optimizers (From top to bottom is SGD, SGDM, Adam, PID, LPF-SGD, HPF-SGD and FuzzyPID).} 
    \label{Figure9} 
\end{figure*}

Though we did not simulate the system response of DDPM, in order to validate that the optimization method should match the learning system, we experiment SGD, SGDM, PID, FuzzyPID, Adam, AdamW \cite{loshchilov2017decoupled}, RAdam \cite{liu2019variance}, NAdam \cite{dozat2016incorporating}, LPF-SGD, and HPF-SGD on the DDPM. The DDPM model only outputs readable samples when using Adam-family optimization methods. GAN also shows the similar preference in the Adam optimizer family. For generating samples from random noise, the Adam-family optimizer outperforms other optimization methods. We speculate on the reason why generation models prefer Adam, and the clue is that Adam has a stronger vibration that can rapidly adjust the learning progress. Other optimization methods, such as, SGDM and PID, have a lower vibration that cannot vigorously adjust the learning and leads to an accumulation of divergence.

\section{Discussion}

Because there is a high consistency between the system responses and the real generation results, we can firstly analyse and simulate the system response of GMs to ensure the stability and convergence, prior to training the clean data on designed optimization methods. Then, GMs can be trained on the clean data using the previously designed reliable component to avoid hallucination on required tasks.

\subsection{The Consistency Between System Response and Generated Samples}

The real learning progress is consistent with the simulation result. GAN satisfies the Adam optimizer, and this is also consistent with the simulation result. The adaptive part of Adam can rapidly adjust the learning process, but other optimizers cannot converge at the end. For particular needs (e.g., Image-to-Image Translation), CycelGAN, this advanced generation system was proposed to generate samples from one data set and to improve its domain adaption on the target data set. Coincidentally, we found that CycleGAN has a preference for PID and FuzzyPID optimizers. The system function of CycleGAN can get a relatively stable solution when using PID and FuzzyPID optimization methods.

\subsection{Laplace Transform Can Help}

Based on control theory, using the Laplace transform can provide with a systematical analysis on GMs, and GMs can tactfully avoid the effect of confident error by utilising a proper optimization component. On the one hand, the system response of GAN cannot converge using SGD, SGDM and PID optimizers. One possible solution is to add an adaptive filter to against the random noise. Thus, Adam outperforms other optimizers on generating samples when using a GAN architecture. Developing a task-matched optimizer according to the system response determines the final performance of GMs. On the other hand, to satisfy various task requirements, learning systems also should become more stable and faster. Classical GAN has a single generator, and the Laplace transform of this single generator has two solutions, which cause to the instability of classical GAN when using other optimization methods. CycleGAN has two generators, which aims to offset the side effect of using single generator. It also considered the utilisation of the cycle loss that ensures the stability of image translation, because the double-generator design
and the cycle loss function can avoid the root solution when using a single generator.

\section{Conclusion}

In this study, we conducted a comprehensive empirical study investigating the connection between control theory and three generation models. By analyzing the system response of these GMs, we explained the rationale behind choosing appropriate optimizers to avoid the hallucination of GMs. Moreover, designing a better learning system under the proper optimization method can also avoid the hallucination of GMs. In our future work, we intend to delve into the control theory and Laplace transform for other ANNs, as well as the development of other optimization components, as we believe the principles of control systems can guide the improvements of all ANNs and their relevant optimization components.

\bibliographystyle{IEEEtran}
\bibliography{ref}

\end{document}